\documentclass[12pt,reqno]{amsart}

\usepackage{fullpage,amsmath,amssymb,amsfonts,mathrsfs,bm,graphicx,hyperref}

\def\bbone{{\mathchoice {\rm 1\mskip-4mu l} {\rm 1\mskip-4mu l}
{\rm 1\mskip-4.5mu l} {\rm 1\mskip-5mu l}}}

\newtheorem{theorem}{Theorem}[section]

\newtheorem{lemma}{Lemma}[section]

\newtheorem{conjecture}{Conjecture}[section]

\begin{document}

\author{Abdelmalek Abdesselam}
\address{Abdelmalek Abdesselam, Department of Mathematics,
P. O. Box 400137,
University of Virginia,
Charlottesville, VA 22904-4137, USA}
\email{malek@virginia.edu}

\title{Log-concavity with respect to the number of orbits for infinite tuples of commuting permutations}

\begin{abstract}
Let $A(p,n,k)$ be the number of $p$-tuples of commuting permutations of $n$ elements whose permutation action results in exactly $k$ orbits or connected components. We formulate the conjecture that, for every fixed $p$ and $n$, the $A(p,n,k)$ form a log-concave sequence with respect to $k$. For $p=1$ this is a well known property of unsigned Stirling numbers of the first kind. As the $p=2$ case, our conjecture includes a previous one by Heim and Neuhauser, which strengthens
a unimodality conjecture for the Nekrasov-Okounkov hook length polynomials. In this article, we prove the $p=\infty$ case of our conjecture. We start from an expression for the $A(p,n,k)$ which follows from an identity by Bryan and Fulman, obtained in the their study of orbifold higher equivariant Euler characteristics. We then derive the $p\rightarrow\infty$ asymptotics. The last step essentially amounts to the log-concavity in $k$ of a generalized Tur\'an number, namely, the maximum product of $k$ positive integers whose sum is $n$. 
\end{abstract}

\maketitle


\section{Introduction and main results}

Let $\langle\sigma_1,\ldots,\sigma_p\rangle$ denote the subgroup of the symmetric group $\mathfrak{S}_n$ generated by the permuations $\sigma_1,\ldots,\sigma_p$. It has an obvious permutation action on the finite set $[n]:=\{1,\ldots,n\}$, with a number of orbits which we will denote by $\kappa(\sigma_1,\ldots,\sigma_p)$. In this article, we will study the numbers $A(p,n,k)$ which count ordered $p$-tuples of permutations $\sigma_1,\ldots,\sigma_p$ of $[n]$, which commute pairwise, and are such that $\kappa(\sigma_1,\ldots,\sigma_p)=k$. We, in particular, propose the following log-concavity conjecture for these numbers.

\begin{conjecture}
For all $p\ge 1$, and for all $n,k$ such that $2\le k\le n-1$, we have the inequality
\[
A(p,n,k)^2\ge A(p,n,k-1)\ A(p,n,k+1)\ .
\]
\label{mainconj}
\end{conjecture}

The $p=1$ case is a classic theorem rather than a conjecture, but we include it in the general statement for the sake of coherence. Indeed, $A(1,n,k)=c(n,k)=|s(n,k)|$, i.e., the unsigned Stirling numbers of the first kind. Their log-concavity property is well known (see, e.g.,~\cite{Abdesselam} and references therein).
The connection of the $p=2$ case to a conjecture by Heim an Neuhauser~\cite[Challenge 3]{HeimN}, related to the Nekrasov-Okounkov formula~\cite{NekrasovO,Westbury}, is explained in the companion article~\cite[\S3]{AbdesselamBDV}.
According to numerical computations by Heim and Neuhauser~\cite{HeimN}, Conjecture \ref{mainconj} is true for $p=2$ and all $n\le 1500$.
In this article, we examine the opposite extreme where $p$ is large, and establish the following result, which one may think of as the ``$p=\infty$'' case of the conjecture.

\begin{theorem}
For all $n\ge 3$, and all $k$ such that $2\le k\le n-1$, we have, as an inequality in the extended real half line $[0,\infty]$,
\begin{equation}
\liminf\limits_{p\rightarrow\infty}
\frac{A(p,n,k)^2}{A(p,n,k-1)\ A(p,n,k+1)}\ \ge 1\ .
\label{liminfeq}
\end{equation}
\label{mainthm}
\end{theorem}

We first note that for all $p\ge 1$ and $n,k$ with $1\le k\le n$, we have $A(p,n,k)\ge1$.
Indeed, one can take $\sigma_1,\ldots,\sigma_p$ all equal to the same permutation given by the cycle $(1, 2, \cdots, n-k+1)$, with one orbit of size $n-k+1$, and $k-1$ other orbits which are singletons. The ratio considered in (\ref{liminfeq}) is well defined and the theorem would follow immediately, upon taking $p\rightarrow\infty$, if one assumes the conjecture is true.
By proving Theorem \ref{mainthm} unconditionally, we thus provide some indirect evidence for Conjecture \ref{mainconj}. In fact, our proof shows the limit and not just the limit inferior exists in $[0,\infty]$, but we do not emphasize this point, which we think is less important than proving the $\ge 1$ lower bound.
The starting point for the proof of the above main theorem is the following precise asymptotic for $A(p,n,k)$ when $p\rightarrow\infty$. Following standard practice, whenever we write $f(p)\sim g(p)$ for positive functions $f,g$ of positive integer variables, we mean the precise statement that 
\[
\lim\limits_{p\rightarrow\infty}\frac{f(p)}{g(p)}=1\ .
\]
In order to state our asymptotic result, we need to define the following quantities.
For $1\le k\le n$, we let
\[
E(n,k):=\max a_1\cdots a_k\ ,
\]
where the maximum is over tuples of integers $a_1,\ldots,a_k\ge 1$ such that $a_1+\cdots+a_k=n$.
It is easy to show that
\begin{equation}
E(n,k)=b^{k-r}(b+1)^r\ ,
\label{Eformula}
\end{equation}
where $n=bk+r$ is the Euclidean division of $n$ by $k$, i.e., $b=\left\lfloor
\frac{n}{k}
\right\rfloor$ is the quotient,
and $r$, $0\le r<k$, is the remainder. Indeed, if for $i\neq j$, we have $a_i\ge a_j+2$, then
\[
(a_i-1)(a_j+1)-a_i a_j=a_i-a_j-1>0\ ,
\] 
so $a_1,\ldots,a_k$ cannot be a maximizing configuration. Such a maximizing configuration would satisfy $|a_i-a_j|\le 1$ for all $i,j$. Namely, there would be some integers $b,r$ with $0\le r<k$ such that $a_i=b$ for $k-r$ values of $i$, and $a_i=b+1$ for the remaining $r$ values of $i$. Since the $a_i$ add up to $n$, the integers $b,r$ are the ones coming from the Euclidean division mentioned above.
This elementary argument also shows that a maximizing configuration is unique, up to permutation. We also note the trivial property $E(n,k)\ge 1$.

We will also need to introduce the multiplicative arithmetic function $H(n)$, for $n\ge 1$, given by
\[
H(n)=\prod_{i=1}^{\ell}\left[
\frac{q_i^{\binom{m_i}{2}-m_i}}{(q_i-1)(q_i^2-1)\cdots(q_i^{m_i}-1)}
\right]\ ,
\]
if the prime factorization of $n$ is $n=q_1^{m_1}\cdots q_{\ell}^{m_{\ell}}$.
For $1\le k\le n$, we now can define
\begin{equation}
F(n,k):=\frac{n!}{r!(k-r)!}\times H(b)^{k-r}
\times H(b+1)^r\ ,
\label{Fdefeq}
\end{equation}
where $b,r$ are the quotient and remainder of the division of $n$ by $k$, as above.

\begin{theorem}
For all $n,k$ with $1\le k\le n$, we have the $p\rightarrow \infty$ asymptotics
\[
A(p,n,k)\sim
F(n,k) \times E(n,k)^p\ .
\]
\label{asympthm}
\end{theorem}

From the exponential behavior in $p$, we therefore see that the log-concavity
with respect to $k$ of the $A(p,n,k)$ is tied to that of the
$E(n,k)$. In fact, we need a bit more, in order to establish our main result, Theorem \ref{mainthm}. This is collected in the next log-concavity result which may be of independent interest, e.g., in the area of extremal combinatorics.

\begin{theorem}
For all $n\ge 3$, and all $k$ such that $2\le k\le n-1$, we have
\[
E(n,k)^2\ge E(n,k-1)\ E(n,k+1)\ .
\]
If, in addition, we have $n\ge k^2-2k+2$, then the strict inequality
\[
E(n,k)^2> E(n,k-1)\ E(n,k+1)
\]
holds.
\label{Ethm}
\end{theorem}

The quantity $E(n,k)$ is an example of generalized Tur\'an number. In the notation of~\cite{AlonS}, we have $E(n,k)=ex(n,K_k,K_{k+1})$, namely the maximal number of $k$-cliques a graph on $n$ vertices can have, if the graph contains no $(k+1)$-clique.
This is indeed a classic result in extremal graph theory resulting from work initiated by Tur\'an~\cite{Turan,Zykov,Erdos}.
By taking graph complements, this is also related to the number of maximal independent sets in a graph, see~\cite{MohrR}.
Our results provide a connection between such numbers, which are of interest in extremal graph theory, and properties of Abelian subgroups of the symmetric group.
Something very similar was already noted in~\cite{BercovM}.
There are works devoted to the study of log-concavity phenomena in the theory of finite Abelian groups, e.g.,~\cite{Butler}, but we are not aware of Conjecture \ref{mainconj}, for general $p$, being stated previously in the literature. 
We also did not see Theorem \ref{Ethm} in previous work, and it raises the question: for what other notable numbers from extremal combinatorics, such as generalized Tur\'an numbers, does one observe similar log-concavity phenomena.

\section{Proofs}

\subsection{Proof of Theorem \ref{asympthm}}
For $p,n\ge 1$, let us define the quantity
\[
B(p,n):=\sum_{s_1|s_2|\cdots|s_{p-1}|n}s_1 s_2\cdots s_{p-1}\ ,
\]
where the summation is over tuples of integers $s_1,\ldots,s_{p-1}\ge 1$ which form an ``arithmetic flag'', i.e., such that $s_1$ divides $s_2$, $s_2$ divides $s_3$,\ldots,
$s_{p-1}$ divides $n$. It is easy to show that $B(p,\cdot)$ is a multiplicative function in the sense of number theory, namely,
\[
B(p,ab)=B(p,a)B(p,b)\ ,
\]
whenever $a,b$ are coprime. The prime power case $n=q^m$, where $q$ is a prime number and $m$ is a nonnegative integer, can be computed explicitly (see~\cite{AbdesselamBDV}), with the result
\begin{equation}
B(p,q^m)=\frac{(q^p-1)(q^{p+1}-1)\cdots(q^{p+m-1}-1)}{(q-1)(q^2-1)\cdots(q^m-1)}\ .
\label{Bformula}
\end{equation}
We now have the following explicit formula for the $A(p,n,k)$ in terms of the $B$ function:
\begin{equation}
A(p,n,k)=\frac{n!}{k!}\times
\sum_{n_1,\ldots,n_k\ge 1}
\bbone\{n_1+\cdots+n_k=n\}\times \prod_{i=1}^{k}
\left[\frac{B(p,n_i)}{n_i}\right]\ ,
\label{Aformula}
\end{equation}
where $\bbone\{\cdots\}$ is the indicator function of the condition between braces.
We refer to~\cite{AbdesselamBDV} for a new proof of this formula, in the spirit of bijective enumerative combinatorics, as well as for an alternate derivation, as a consequence of a previous generating function identity by Bryan and Fulman~\cite{BryanF}, which they discovered in the context of orbifold higher equivariant Euler characteristics.
The number $B(p,n)$ is a well studied object (see, e.g.,~\cite[Ch. 15]{LubotzkyS}, as well as~\cite{Solomon}) since it counts subgroups of $\mathbb{Z}^p$ which are of index $n$.
Our derivation of the $p\rightarrow\infty$ asymptotics for $A(p,n,k)$ is a rather straightforward consequence of (\ref{Bformula}) and (\ref{Aformula}).
In the prime power case, we get
\[
B(p,q^m)\sim 
\frac{q^{\binom{m}{2}}}{(q-1)(q^2-1)\cdots(q^{m}-1)}\times
q^{mp}\ ,
\]
when $p\rightarrow\infty$. 
As a result, for $n$ with prime factorization $n=q_1^{m_1}\cdots q_{\ell}^{m_{\ell}}$,
we get $\frac{B(p,n)}{n}\sim H(n)\times n^p$.
A term in the sum (\ref{Aformula}) is asymptotically equivalent to
\[
C\times (n_1\cdots n_k)^p\ ,
\]
where $C$ is a constant depending on $n_1,\ldots,n_k$, but not on $p$.
Clearly, the asymptotic behavior of the sum is dictated by the dominant terms with the highest rate of exponential growth, i.e., the terms which maximize the product $n_1\cdots n_k$.
By the discussion of $E(n,k)$ in the introduction, there are $\binom{k}{r}$ such terms where $k-r$ of the $n_i$ are equal to $b$, and $r$ of the $n_i$ are equal to $b+1$. Again, the notation is $b=\left\lfloor\frac{n}{k}\right\rfloor$, and $r=n-bk$.
After making the constant $C$ explicit, for these dominant terms, and cleaning up, we arrive at the asymptotics in Theorem \ref{asympthm}.
\qed

\subsection{Proof of Theorem \ref{Ethm}}
The proof uses two related lemmas, the second being a reinforcement of the first, for a smaller range of parameters.

\begin{lemma}
For all $n,k$ with $2\le k\le n$, we have
\[
\frac{E(n,k)}{E(n,k-1)}\le 
\frac{E(n+1,k)}{E(n+1,k-1)}\ .
\]
In other words, for all fixed $k\ge 2$, the expression $\frac{E(n,k)}{E(n,k-1)}$ is a nondecreasing function of $n$, in the range $n\ge k$.
\label{weaklem}
\end{lemma}

\noindent{\bf Proof:}
Pick a maximizing configuration $a_1,\ldots,a_k\ge 1$, such that
$E(n,k)=a_1\cdots a_k$, with $a_1+\cdots+a_k=n$.
Also pick a maximizing configuration $b_1,\ldots,b_{k-1}\ge 1$, such that
$E(n+1,k-1)=b_1\cdots b_{k-1}$, with $b_1+\cdots+b_{k-1}=n+1$.
We argue by contradiction, in order to show $\exists i\in [k]$ and $\exists j\in[k-1]$, such that $a_i+1\le b_j$. Indeed, we would otherwise have $a_i\ge b_j$ for all $i,j$.
This would imply
\[
n=a_1+\cdots+a_k\ge b_1+\cdots+b_{k-1}+a_k\ge n+2\ ,
\]
which is absurd. Without loss of generality, one may renumber the $a$'s and $b$'s so that the produced $i,j$ are both equal to $1$, and therefore $a_1+1\le b_1$.
Then, we have
\[
(a_1+1)(b_1-1)\ge a_1 b_1\ ,
\]
which implies
\[
(a_1+1)a_2\cdots a_k\times (b_1-1)b_2\cdots b_{k-1}\ge
a_1 a_2\cdots a_k\times b_1 b_2\cdots b_{k-1}\ .
\]
Since, by definition, $E(n+1,k)\ge (a_1+1)a_2\cdots a_k$,
and $E(n,k-1)\ge (b_1-1)b_2\cdots b_{k-1}$, we immediately arrive at
\begin{equation}
E(n+1,k)E(n,k-1)\ge E(n,k)E(n+1,k-1)\ ,
\label{finalineq}
\end{equation}
and the lemma follows.
\qed

\begin{lemma}
For all $k\ge 2$, and all $n$ such that $n\ge k^2-2k+1$, we have
\[
\frac{E(n,k)}{E(n,k-1)}<
\frac{E(n+1,k)}{E(n+1,k-1)}\ .
\]
\label{strictlem}
\end{lemma}

\noindent{\bf Proof:}
The proof is a repeat of the previous one, but where we now make the stronger claim $\exists i\in [k]$ and $\exists j\in[k-1]$, such that $a_i+2\le b_j$. Following the same steps as before, we would get $(a_1+1)(b_1-1)> a_1 b_1$ and the final inequality (\ref{finalineq}) would become strict.
To see why the new claim holds, assume for the sake of contradiction that $a_i\ge b_j-1$, for all $i\in[k]$ and $j\in[k-1]$.
We then have
\[
n=a_1+\cdots+a_k\ge k\times\min a_i
\ge k(\max b_j)-k\ .
\]
On the other hand,
\[
n+1=b_1+\cdots+b_{k-1}\le (k-1)\times \max b_j\ .
\]
From the previous inequalities, one would deduce
\[
n\ge k\times \frac{n+1}{k-1}-k\ ,
\]
i.e., $n\le k^2-2k$, but this is ruled out by the hypothesis $n\ge k^2-2k+1$. Hence our claim holds, and the lemma is proved.
\qed

We now address the proof of the first part of Theorem \ref{Ethm} and consider $n,k$ with $2\le k\le n-1$.
We pick a maximizing configuration $a_1,\ldots,a_{k-1}\ge 1$, such that
$E(n,k-1)=a_1\cdots a_{k-1}$, with $a_1+\cdots+a_{k-1}=n$.
We also pick a maximizing configuration $b_1,\ldots,b_{k+1}\ge 1$, such that
$E(n,k+1)=b_1\cdots b_{k+1}$, with $b_1+\cdots+b_{k+1}=n$.
Since
\[
a_1+\cdots+a_{k-1}=n>b_1+\cdots+b_{k}\ge k\ge 2\ ,
\]
Lemma \ref{weaklem} implies
\begin{equation}
\frac{E(n,k)}{E(n,k-1)}\ge\frac{E(b_1+\cdots+b_k,k)}{E(b_1+\cdots+b_k,k-1)}\ .
\label{intermineq}
\end{equation}
Pick a maximizing configuration $c_1,\ldots,c_k\ge 1$, such that
$E(n,k)=c_1\cdots c_k$, with $c_1+\cdots+c_k=n$.
Also pick a maximizing configuration $d_1,\ldots,d_{k-1}\ge 1$, such that
$E(b_1+\cdots+b_k,k-1)=d_1\cdots d_{k-1}$, with $d_1+\cdots+d_{k-1}=b_1+\cdots+b_k$.
We note that $E(b_1+\cdots+b_k,k)=b_1\cdots b_k$, otherwise the larger collection $b_1,\ldots,b_{k+1}$ would not be a maximizing configuration for the computation of $E(n,k+1)$.
From (\ref{intermineq}), we obtain
\[
c_1\cdots c_k d_1\cdots d_{k-1}\ge a_1\cdots a_{k-1} b_1\cdots b_{k}\ .
\]
After multiplying by $b_{k+1}$, we arrive at
\begin{equation}
\left(c_1\cdots c_k\right)\times
\left(d_1\cdots d_{k-1} b_{k+1}\right)
\ge
\left(a_1\cdots a_{k-1}\right)\times
\left(b_1\cdots b_{k+1}\right)\ .
\label{switchineq}
\end{equation}
Since, by construction, $E(n,k)=c_1\cdots c_k$, and, by definition, $E(n,k)\ge d_1\cdots d_{k-1} b_{k+1}$, the inequality (\ref{switchineq}) yields
\[
E(n,k)^2\ge E(n,k-1) E(n,k+1)\ .
\]
For the second part of the theorem, with the strict inequality, we see that, in the previous argument, we only need to make sure that the inequality in (\ref{intermineq}) is strict.
This follows from Lemma \ref{strictlem} because $n-1\ge k^2-2k+1$, by assumption on $n$.
\qed

\subsection{Proof of Theorem \ref{mainthm}}

We now pick $n\ge 3$ and $k$ such that $2\le k\le n-1$, and perform the following three Euclidean divisions
\begin{eqnarray}
n & = & (k-1)b_1 +r_1\nonumber\\
n & = & kb_2 +r_2\nonumber\\
n & = & (k+1)b_3 +r_3\ .
\label{Eucldiv}
\end{eqnarray}
Clearly, the integers $b_1,b_2,b_3$ satisfy $b_1\ge b_2\ge b_3\ge 1$, and by definition, we have $0\le r_1<k-1$, $0\le r_2<k$, and $0\le r_3<k+1$.
We form the following ratios
\begin{eqnarray*}
R_E(n,k) & := & \frac{E(n,k)^2}{E(n,k-1)\ E(n,k+1)} \\
R_F(n,k) & := &  \frac{F(n,k)^2}{F(n,k-1)\ F(n,k+1)} \\
R_{F,1}(n,k) & := & \frac{r_1!\ (k-1-r_1)!\ r_3!\ (k+1-r_3)!}{r_2!^2\ (k-r_2)!^2} \\
R_{F,2}(n,k) & := & \frac{H(b_2)^{2(k-r_2)}\ H(b_2+1)^{2r_2}}{H(b_1)^{k-1-r_1}\ 
H(b_1+1)^{r_1}\ H(b_3)^{k+1-r_3}\ H(b_3+1)^{r_3}}\ . 
\end{eqnarray*}
From (\ref{Fdefeq}), we see that $R_F(n,k)=R_{F,1}(n,k)R_{F,2}(n,k)$,
and Theorem \ref{asympthm} entails
\[
\frac{A(p,n,k)^2}{A(p,n,k-1)\ A(p,n,k+1)}\sim
R_F(n,k)\times R_E(n,k)^p
\ ,
\]
when $p\rightarrow\infty$.
Therefore, Theorem \ref{mainthm} amounts to showing either $R_E(n,k)>1$, or
$R_E(n,k)=1$ together with $R_{F}(n,k)\ge 1$. By Theorem \ref{Ethm}, we already know $R_E(n,k)\ge 1$ always holds. We also have $R_E(n,k)>1$, for $k$ less than roughly $\sqrt{n}$. For the opposite regime, we need the following lemma, in order to control the drops in the sequence of quotients $b_1,b_2,b_3$.

\begin{lemma}
If $n\le k^2-k$, then $b_1-b_2\le 1$ and $b_2-b_3\le 1$.
\label{droplem}
\end{lemma}

\noindent{\bf Proof:}
We have
\[
b_1-b_2=\left\lfloor\frac{n}{k-1}\right\rfloor
-\left\lfloor\frac{n}{k}\right\rfloor
<\left(\frac{n}{k-1}\right)-\left(\frac{n}{k}-1\right)
=\frac{n}{k(k-1)}+1\le 2\ .
\]
Since one of the inequalities in the above chain is strict, and we are bounding integers, we, in fact, have $b_1-b_2\le 1$.
By the same argument, with $k+1$ instead of $k$, we get $b_2-b_3\le 1$,
since the hypothesis also implies $n\le k^2+k$.
\qed

In order to compare some products of factorials, we will also need the next easy lemma.

\begin{lemma}
Let $\alpha,\beta,\gamma,\delta$ be nonnegative numbers with $\alpha+\beta=\gamma+\delta$. If $\min(\alpha,\beta)\ge\min(\gamma,\delta)$, then $\alpha!\beta!\le\gamma!\delta!$.
\label{binomlem}
\end{lemma}

\noindent{\bf Proof:}
Let $S:=\alpha+\beta=\gamma+\delta$. As is well known (see, e.g.,~\cite{Zeilberger}), the binomial coefficients $\binom{S}{x}$ are nondecreasing in the range $0\le x\le \left\lfloor\frac{S}{2}\right\rfloor$.
Without loss of generality, we may assume $\alpha\le \beta$ and $\gamma\le \delta$, and therefore $0\le\gamma\le\alpha\le\frac{S}{2}$.
The lemma then follows from the inequality $\binom{S}{\gamma}\le\binom{S}{\alpha}$.
\qed

We now proceed with the proof of Theorem \ref{mainthm}, and break it into six cases.

\medskip\noindent
{\bf 1st case:} When $n\ge k^2-2k+2$.
Then, the second part of Theorem \ref{Ethm} gives $R_E(n,k)>1$ and we are done.

\medskip
In the remaining five cases, we will assume $n\le k^2-2k+1$. Since $k^2-2k+1\le k^2-k$, Lemma \ref{droplem} applies and forces the drops between successive $b_i$'s to be at most one. This gives four possibilities for the triple $(b_1,b_2,b_3)$, one of which leading to two subcases.

\medskip\noindent
{\bf 2nd case:} When $(b_1,b_2,b_3)=(b,b,b)$ for some $b$.
From the Euclidean divisions (\ref{Eucldiv}), one must have $2r_2=r_1+r_3$.
Then, a quick computation shows that
$R_E(n,k)=1$ and $R_{F,2}(n,k)=1$.
The nontrivial factor is
\[
R_{F,1}(n,k)=\frac{(n-(k-1)b)!\ ((k-1)(b+1)-n)!\ (n-(k+1)b)!\ ((k+1)(b+1)-n)!}{(n-kb)!^2\ (k(b+1)-n)!^2}\ .
\]
Since $n-kb\ge n-(k+1)b$, Lemma \ref{binomlem} implies
\[
(n-kb)!^2\le (n-(k-1)b)!\ (n-(k+1)b)!\ . 
\]
Since $k(b+1)-n\ge (k-1)(b+1)-n$, Lemma \ref{binomlem} implies
\[
(k(b+1)-n)!^2\le ((k-1)(b+1)-n)!\ ((k+1)(b+1)-n)!\ . 
\]
As a result, $R_{F,1}(n,k)\ge 1$, and thus $R_F(n,k)\ge 1$, as wanted.

\medskip\noindent
{\bf 3rd case:} When $(b_1,b_2,b_3)=(b+1,b,b)$ for some $b$, and $r_1>0$.
Using
\begin{eqnarray*}
r_1 & = & n-(b+1)(k-1)\\
r_2 & = & n-bk\\
r_3 & = & n-b(k+1)\ ,
\end{eqnarray*}
and the formula (\ref{Eformula}), we arrive, after simplification, at
\[
R_E(n,k)=\left[\frac{(b+1)^2}{b(b+2)}\right]^{n-(k-1)(b+1)}\ .
\]
Note that $(b+1)^2>b(b+2)$.
Since $r_1>0$, the exponent is strictly positive, which gives the desired conlcusion $R_E(n,k)>1$.

\medskip\noindent
{\bf 4th case:} When $(b_1,b_2,b_3)=(b+1,b,b)$ for some $b$, and $r_1=0$.
The discussion in the previous case shows that we now have $R_E(n,k)=1$.
A straightforward computation, taking into account the relation $n=(k-1)(b+1)$, shows that $R_{F,2}(n,k)=1$ and
\[
R_{F,1}(n,k)=\frac{(k-1)!\ (k-2b-1)!\ (2b+2)!}{(k-b-1)!^2\ (b+1)!^2}\ .
\]
Since $k-b-1\ge k-2b-1$, Lemma \ref{binomlem} implies
\[
(k-b-1)!^2\le (k-1)!\ (k-2b-1)!\ ,
\]
and therefore
\[
R_{F,1}(n,k)\ge \binom{2b+2}{b+1}\ge 1\ .
\]
As a result, we obtain$R_{F}(n,k)\ge 1$, as desired.

\medskip\noindent
{\bf 5th case:} When $(b_1,b_2,b_3)=(b+1,b+1,b)$ for some $b$.
Using
\begin{eqnarray*}
r_1 & = & n-(b+1)(k-1)\\
r_2 & = & n-(b+1)k\\
r_3 & = & n-b(k+1)\ ,
\end{eqnarray*}
and the formula (\ref{Eformula}), we arrive, after simplification, at
\[
R_E(n,k)=\left[\frac{(b+1)^2}{b(b+2)}\right]^{(k+1)(b+1)-n}\ .
\]
We have $(k+1)(b+1)-n=k+1-r_3>0$ by definition of the remainder $r_3$, and therefore,
$R_E(n,k)>1$ as hoped for.

\medskip\noindent
{\bf 6th case:} When $(b_1,b_2,b_3)=(b+2,b+1,b)$ for some $b$.
Using
\begin{eqnarray*}
r_1 & = & n-(b+2)(k-1)\\
r_2 & = & n-(b+1)k\\
r_3 & = & n-b(k+1)\ ,
\end{eqnarray*}
and the formula (\ref{Eformula}), we arrive, after simplification, at
\[
R_E(n,k)=\left[\frac{(b+1)^2}{b(b+2)}\right]^{(k+1)(b+1)-n}\times
\left[\frac{(b+2)^2}{(b+1)(b+3)}\right]^{n-(k-1)(b+2)}
\ .
\]
Both ratios being exponentiated are $>1$.
As for the exponents, we have $(k+1)(b+1)-n=k+1-r_3>0$ and $n-(k-1)(b+2)=r_1\ge 0$, by definition of the remainders $r_1,r_3$.
This again implies $R_E(n,k)>1$ and we are done.

\qed

\bigskip
\noindent{\bf Acknowledgements:}
{\small
The author thanks Ken Ono for introducing him to the Nekrasov-Okounkov formula and the Amdeberhan-Heim-Neuhauser unimodality conjecture.
The author also thanks MathOverflow users Nate, Tony Huynh, and Carlo Beenakker for useful information regarding the quantity $E(n,k)$. The author thanks his mathematical quantum field theory students Pedro Brunialti, Tristan Doan, and Philip Velie for the collaboration~\cite{AbdesselamBDV}. Finally, we thank Bernhard Heim, Markus Neuhauser, and Alexander Thomas for a careful reading of the first version of this article.
}


\begin{thebibliography}{999}

\bibitem{Abdesselam}
A. Abdesselam, A local injective proof of log-concavity for increasing spanning forests.
Discrete Math. {\bf 346} (2023), no. 12, Paper No. 113651.

\bibitem{AbdesselamBDV}
A. Abdesselam, P. Brunialti, T. Doan, and P. Velie, A bijection for tuples of commuting permutations and a log-concavity conjecture. Preprint 2023.

\bibitem{AlonS}
N. Alon, and C. Shikhelman,
Many $T$ copies in $H$-free graphs.
J. Combin. Theory Ser. B {\bf 121} (2016), 146--172.

\bibitem{BercovM}
R. Bercov, and L. Moser,
On Abelian permutation groups.
Canad. Math. Bull. {\bf 8} (1965), 627–630.

\bibitem{BryanF}
J. Bryan, and J. Fulman,
Orbifold Euler characteristics and the number of commuting $m$-tuples in the symmetric groups.
Ann. Comb. {\bf 2} (1998), no. 1, 1--6.

\bibitem{Butler}
L. M. Butler,
A unimodality result in the enumeration of subgroups of a finite abelian group.
Proc. Amer. Math. Soc. {\bf 101} (1987), no. 4, 771--775.

\bibitem{Erdos}
P. Erd\H{o}s,
On the number of complete subgraphs contained in certain graphs.
Magyar Tud. Akad. Mat. Kutat\'{o} Int. K\"{o}zl. {\bf 7} (1962), 459–464.

\bibitem{HeimN}
B. Heim, and M. Neuhauser,
Horizontal and vertical log-concavity.
Res. Number Theory {\bf 7} (2021), no. 1, Paper No. 18, 12 pp.

\bibitem{LubotzkyS}
A. Lubotzky, and D. Segal,
{\it Subgroup Growth}.
Progr. Math. {\bf 212},
Birkh\"{a}user Verlag, Basel, 2003. 

\bibitem{MohrR}
E. Mohr, and D. Rautenbach,
On the maximum number of maximum independent sets.
Graphs Combin. {\bf 34} (2018), no. 6, 1729--1740.

\bibitem{NekrasovO}
N. A. Nekrasov, and A. Okounkov,
Seiberg-Witten theory and random partitions. In: {\it The Unity of Mathematics}. In honor of the ninetieth birthday of I. M. Gelfand. Papers from the conference held in Cambridge, MA, August 31--September 4, 2003, pp. 525--596. Eds: P. Etingof, V. Retakh, and I. M. Singer,
Progr. Math. {\bf 244},
Birkh\"{a}user Boston, Inc., Boston, MA, 2006.

\bibitem{Solomon}
L. Solomon,
Partially ordered sets with colors.
In: {\it Relations Between Combinatorics and Other Parts of Mathematics}.
Proceedings of the Symposium in Pure Mathematics of the American Mathematical Society held at the Ohio State University, Columbus, Ohio, March 20--23, 1978. 
pp. 309--329.
Ed: D. K. Ray-Chaudhuri,
Proc. Sympos. Pure Math. {\bf XXXIV},
American Mathematical Society, Providence, RI, 1979.

\bibitem{Turan}
P. Tur\'an,
Eine Extremalaufgabe aus der Graphentheorie.
Mat. Fiz. Lapok {\bf 48} (1941), 436--452.

\bibitem{Westbury}
B. W. Westbury,
Universal characters from the Macdonald identities.
Adv. Math. {\bf 202} (2006), no. 1, 50--63.

\bibitem{Zeilberger}
D. Zeilberger,
$\binom{5}{2}$ proofs that $\binom{n}{k} \le \binom{n}{k+1}$ if $k<n/2$.
Preprint  arXiv:1003.1273[math.CO], 2010.

\bibitem{Zykov}
A. A. Zykov, On some properties of linear complexes.
Amer. Math. Soc. Translation {\bf 1952} (1952), no. 79, 33 pp.
Translated from Mat. Sbornik N.S. {\bf 24/66} (1949), 163--188.

\end{thebibliography}
\end{document}